\documentclass[11pt]{article}

\usepackage{amsmath,amssymb}
\usepackage{graphicx}  
\usepackage{pict2e}
\usepackage{subcaption}
\usepackage{algorithm}
\usepackage[noend]{algpseudocode}
\usepackage{tikz}
\usetikzlibrary{math}
\usepackage{listofitems}
\usepackage{xcolor}

\newcommand{\sfrac}[2]{{\textstyle\frac{#1}{#2}}}

\newcommand{\eps}{\varepsilon}
\newcommand{\bP}{\mathbf P}
\newcommand{\bx}{\mathbf x}

\newcommand{\comment}[1]{}

\begin{document}

\title{A Real-World Markov Chain arising in Recreational Volleyball}
 \author{David J. Aldous\thanks{Department of Statistics,
 367 Evans Hall \#\  3860,
 U.C. Berkeley CA 94720;  aldous@stat.berkeley.edu;
  www.stat.berkeley.edu/users/aldous.}
\and Madelyn Cruz\thanks{University of the Philippines - Diliman; yumadelynesther@gmail.com}
}

 \maketitle
 
  \begin{abstract}
  Card shuffling models have provided simple motivating examples for the mathematical theory of mixing times for Markov chains.
  As a complement, we introduce a more intricate realistic model of a certain observable real-world scheme for mixing human players onto teams.
  We quantify numerically the effectiveness of this mixing scheme over the 7 or 8 steps performed in practice.
  We give a combinatorial proof of the non-trivial fact that the chain is indeed irreducible.
  \end{abstract}
  
  {\bf Key words:} Markov chain, mixing time.
  
  {\bf MSC subject classification:} 60J10
 
\section{Introduction}
In introducing Markov chains at some elementary level, the first author always found it difficult to give motivating examples with a real-world story, a plausible probability model, and a fairly rich mathematical structure.
Then he realized that he was a regular participant in one such story. 
A first thought was to write out the model for possible use as an instructional example in an introductory lecture. As often happens,  things turned out  
to be more complicated than first imagined, so it was re-purposed as a basis for a challenging undergraduate project to study further aspects of the model. 
The second author took up the challenge.
Some remaining questions that could be used for undergraduate projects are mentioned in section \ref{sec:projects}.

\section{The model}

The story concerns recreational volleyball, in a ``drop-in" setting without fixed teams, and where one wants the team compositions to change from game to game, both as socialization and to avoid persistent large differences in team skill levels. 
Specifically, there are 24 people, and at each stage, there are two ongoing games on two courts, each game between two teams, each team with 6 players on a half-court.
Over the 2 hour period there will be 7 or 8 successive such stages, everyone always playing.
The rule\footnote{Actually used in the gym where the first author plays; I don't know how common it is.}  for changing team composition is very simple, exploiting a particular incidental feature
of volleyball:

{\em At the end of one stage, the players in the back row of each team stay in these positions
for the start of the next game, while the front row players move (clockwise in the gym) to the same positions in the next quadrant.}

See Figure \ref{Fig:1}.
\setlength{\unitlength}{0.11in}
\begin{figure}[h]
\begin{center}
\begin{picture}(18,36)
\put(0.3,0){\line(1,0){6}}
\put(0.3,4.3){\line(1,0){6}}
\put(0.3,8.6){\line(1,0){6}}
\put(0.3,0){\line(0,1){8.6}}
\put(6.3,0){\line(0,1){8.6}}
\put(1,1){c}
\put(3,1){d}
\put(5,1){e}
\put(1,3){$\beta$}
\put(3,3){$\alpha$}
\put(5,3){$\phi$}
\put(1,5){f}
\put(3,5){a}
\put(5,5){b}
\put(1,7){F}
\put(3,7){E}
\put(5,7){D}
\put(21,5){start game 1}
\put(12.3,0){\line(1,0){6}}
\put(12.3,4.3){\line(1,0){6}}
\put(12.3,8.6){\line(1,0){6}}
\put(12.3,0){\line(0,1){8.6}}
\put(18.3,0){\line(0,1){8.6}}
\put(13,1){$\gamma$}
\put(15,1){$\delta$}
\put(17,1){$\eps$}
\put(13,3){i}
\put(15,3){vi}
\put(17,3){v}
\put(13,5){A}
\put(15,5){B}
\put(17,5){C}
\put(13,7){iv}
\put(15,7){iii}
\put(17,7){ii}
\put(0.3,15){\line(1,0){6}}
\put(0.3,19.3){\line(1,0){6}}
\put(0.3,23.6){\line(1,0){6}}
\put(0.3,15){\line(0,1){8.6}}
\put(6.3,15){\line(0,1){8.6}}
\put(1,16){c}
\put(3,16){d}
\put(5,16){e}
\put(1,18){b}
\put(3,18){a}
\put(5,18){f}
\put(1,20){A}
\put(3,20){B}
\put(5,20){C}
\put(1,22){F}
\put(3,22){E}
\put(5,22){D}
\put(21,20){end game 0}
\put(12.3,15){\line(1,0){6}}
\put(12.3,19.3){\line(1,0){6}}
\put(12.3,23.6){\line(1,0){6}}
\put(12.3,15,0){\line(0,1){8.6}}
\put(18.3,15,0){\line(0,1){8.6}}
\put(13,16){$\gamma$}
\put(15,16){$\delta$}
\put(17,16){$\eps$}
\put(13,18){$\beta$}
\put(15,18){$\alpha$}
\put(17,18){$\phi$}
\put(13,20){v}
\put(15,20){vi}
\put(17,20){i}
\put(13,22){iv}
\put(15,22){iii}
\put(17,22){ii}
\put(0.3,30){\line(1,0){6}}
\put(0.3,34.3){\line(1,0){6}}
\put(0.3,38.6){\line(1,0){6}}
\put(0.3,30){\line(0,1){8.6}}
\put(6.3,30){\line(0,1){8.6}}
\put(1,31){e}
\put(3,31){f}
\put(5,31){a}
\put(1,33){d}
\put(3,33){c}
\put(5,33){b}
\put(1,35){B}
\put(3,35){C}
\put(5,35){D}
\put(1,37){A}
\put(3,37){F}
\put(5,37){E}
\put(21,35){start game 0}
\put(12.3,30){\line(1,0){6}}
\put(12.3,34.3){\line(1,0){6}}
\put(12.3,38.6){\line(1,0){6}}
\put(12.3,30,0){\line(0,1){8.6}}
\put(18.3,30,0){\line(0,1){8.6}}
\put(13,31){$\eps$}
\put(15,31){$\phi$}
\put(17,31){$\alpha$}
\put(13,33){$\delta$}
\put(15,33){$\gamma$}
\put(17,33){$\beta$}
\put(13,35){ii}
\put(15,35){iii}
\put(17,35){iv}
\put(13,37){i}
\put(15,37){vi}
\put(17,37){v}
\end{picture}
\caption{One step of the big chain.  As in football, {\em positions} are relative to the way a team is facing:  at start of game 0 both $b$ and $B$ are {\em front right} in their teams.}
\label{Fig:1}
\end{center}
\end{figure}
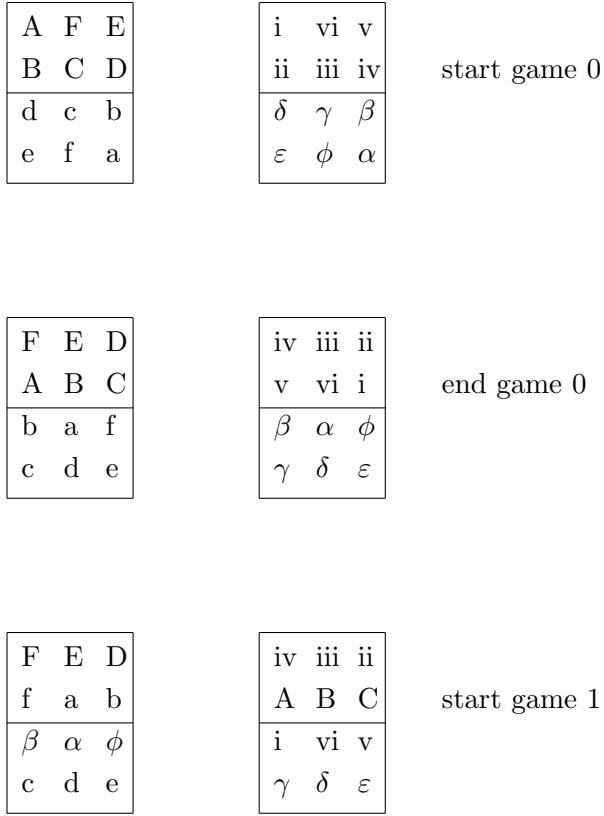
The key point is that in volleyball, there are 6 ``positions"\footnote{By convention numbered 1 to 6, starting in serving position (back right, as facing the net) and ordered counter-clockwise.
Because players rotate clockwise, this indicates serving order.}, and players rotate one position each time their team regains the serve,
and this happens a {\em random} number of times during a game.  So, relative to initial positions, the 3 players who finish in the front row
will in fact be (to a good approximation) a {\em uniform} random choice over the 6 possibilities of 3 adjacent players, and so we model this as a uniform random choice.

To complete a mathematical model, note that the number of one-position rotations of two opposing teams can differ (because they alternate rotations) by at most one.
So,  independently for the two courts, we model the final positions of players in opposing teams in a game as rotations by $(C_1,C_2)$ where $C_1$ is 
uniform\footnote{In a more detailed model, each team wins a random Geometric(1/2) number of points between rotates, and the game ends when one team reaches 25 points.  Then the number of rotations modulo 6 does indeed have distribution close to uniform.} 
on $\{0,1,\ldots,5\}$ and 
$C_2 = C_1 - 1 + \mathrm{Binomial}(2,1/2)$ modulo 6, the Binomial term reflecting randomness of the initial serving team and the final serving team. 

This specifies a ``big" Markov chain on the 24! states (assignment of players to positions).
One step of this {\em big chain} is from the starting positions in one game to  the starting positions in the next game,
as illustrated in  Figure \ref{Fig:1}.
This model is conceptually loosely related to some card shuffling models\footnote{See further discussion in section \ref{sec:final}.} such as 
those in \cite{shuffle,MCMT} in that a ``rotation" 
of team players corresponds to a cut-shuffle of a 6-card deck.
But unlike playing cards, the volleyball players  care about their positions relative to other players for various 
reasons\footnote{Friendly rivalry between spikers/blockers;  more talented setters enable sophisticated fast plays; attractive members of opposite sex; \ldots} and this suggests actual observables for study in the model.

\section{Results}

The central, albeit vague, question is 
\begin{quote} 
how effective is this scheme at mixing up the teams?
\end{quote}
In a lecture course, this would provide a real-world example  for later discussion of the {\em mixing times} topic. 
It seems intuitively obvious that this scheme would mix perfectly in the long run.
What does that mean? 
First observe that the chain is doubly stochastic.  
To see why, consider the step illustrated in Figure \ref{Fig:1}.  
This takes the ``start game 0" configuration $\bx_0$ to a certain ``start game 1" configuration $\bx_1$, via a series of rotations.
By reversing that series, one sees that there exists a configuration $\bx_{-1}$ which, from the same (forwards) series of rotations, takes $\bx_{-1}$ to $\bx_0$.
This leads to (for fixed $\bx_0$) a bijection between possible configurations $\bx_1$ and $\bx_{-1}$ which preserves transition probabilities; 
which in turn implies the doubly stochastic property.
And that property implies that the uniform distribution on all $24!$ states is a stationary distribution for the chain.\footnote{More generally, any card-shuffling scheme in which the shuffling rule depends only on the ranks (positions within deck) and not on the labels of the cards will be doubly stochastic.}

 Basic finite Markov chain theory\footnote{In many textbooks such as \cite{haggstrom,privault}.} 
  identifies ``mix perfectly in the long run" with {\em irreducible and aperiodic}, which implies convergence of time-$t$ distributions to a unique 
  stationary distribution, which in our model must be the uniform distribution on all $24!$ states. 
  Theory also tells us that {\em irreducible} is equivalent to the property
\begin{quote} 
 the directed graph of all possible transitions on the 24! states is strongly connected.
\end{quote}
This property is purely combinatorial -- the numerical values of the non-zero transition probabilities do not matter.  
So our first goal is to prove irreducibility and aperiodicity.
We give a constructive proof of irreducibility in section \ref{sec:irreducible}, and a proof of aperiodicity in section \ref{sec:aperiodic}.
Our proofs are rather complicated, and no doubt there exist simpler proofs.

Our second set of results concern numerical calculation or simulation of statistics relating to the realistic short term in this story -- 7 or 8 steps.
Some basic observables involve the {\em friend chain} indicating the relative positions of two players.
A variety of numerical results are shown in section \ref{sec:friend}. 
For instance, if your friend does not start on your team, then the probability that you are never on the same team over 8 games varies between 
$0.251$ and $0.403$ depending on initial relative positions (Table \ref{tab:notsame}).
A pedagogic point is that, for numerical calculations, 
we don't want to work with $24! \times 24!$ transition matrices, but instead exploit symmetry to reduce to question-specific small-state chains. 
For instance the way a given player moves between games is simple: with chance $1/2$ they stay, with chance $1/2$ they move to the next quadrant. In jargon, the lazy cyclic walk \cite{MCMT}.

The bottom line is that, as regards simple observables, this scheme does a reasonable job of  mixing up the teams over 8 games. 
However, the central point of sophisticated {\em mixing time} theory \cite{MCMT}
is to go beyond the unquantified ``eventually" implied by irreducibility,
and instead to quantify when the step-$t$ distribution is close to (in our case) the uniform distribution.
The usual quantification involves {\em variation distance} between distributions on the $24!$ states,
and this is the context of the famous Bayer-Diaconis result \cite{bayer} for riffle shuffles,
informally called ``7 shuffles suffice" \cite{NYT}.
Studying variation distance for our big chain, either numerically of via analytic bounds, remains a challenging open problem.
We give some preliminary observations in section \ref{sec:mixing}.

\newpage
\section{The big chain is irreducible}
\label{sec:irreducible}
\subsection{Notation}
To prove that the ``big" chain is irreducible, we will show that it is possible to move from any one given state to any other given state 
via some sequence of allowable transitions of the chain.

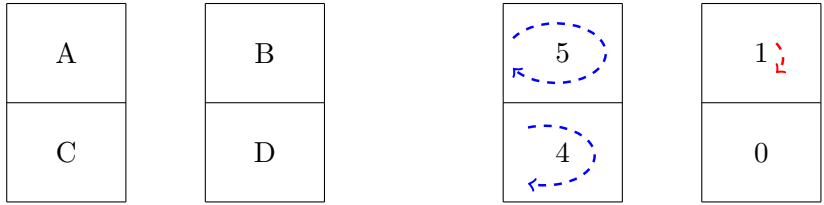
\begin{figure}[h]
    \centering
    \begin{tikzpicture}[scale = 0.66]

    \draw[blue,dashed,line width=.2ex] (0,-0.2) [->] arc (150:-150:1 and 0.6);
    \draw[blue,dashed,line width=.2ex] (0.3,-2) [->] arc (110:-110:1 and 0.6);
     \node at (1,-0.5) {5} ; 
     \node at (1,-2.5) {4} ; 
     \draw (-0.2,-1.5) -- (2.2,-1.5);
     \draw (-0.2,0.5) -- (-0.2,-3.5);
     \draw (-0.2,0.5) -- (2.2,0.5);
     \draw (-0.2,-3.5) -- (2.2,-3.5);
     \draw (2.2,0.5) -- (2.2,-3.5); 
     
     \draw (3.8,-1.5) -- (6.2,-1.5);
     \draw (3.8,0.5) -- (3.8,-3.5);
     \draw (3.8,0.5) -- (6.2,0.5);
     \draw (3.8,-3.5) -- (6.2,-3.5);
     \draw (6.2,0.5) -- (6.2,-3.5); 
     \node at (5,-0.5) {1} ; 
     \node at (5,-2.5) {0} ; 
      
     \draw[red,dashed,line width=.2ex] (5.3,-0.3) [->] arc (30:-30:1 and 0.6);

      \draw (-6.2,-1.5) -- (-3.8,-1.5);  
     \draw (-6.2,0.5) -- (-6.2,-3.5);
     \draw (-6.2,0.5) -- (-3.8,0.5);
     \draw (-6.2,-3.5) -- (-3.8,-3.5);
     \draw (-3.8,0.5) -- (-3.8,-3.5); 
     
      \draw (-10.2,-1.5) -- (-7.8,-1.5);
      \draw (-10.2,0.5) -- (-10.2,-3.5);
     \draw (-10.2,0.5) -- (-7.8,0.5);
     \draw (-10.2,-3.5) -- (-7.8,-3.5);
     \draw (-7.8,0.5) -- (-7.8,-3.5); 
     
      \node at (-9,-0.5) {A} ; 
     \node at (-9,-2.5) {C} ; 
     \node at (-5,-0.5) {B} ; 
     \node at (-5,-2.5) {D} ; 
    
    \end{tikzpicture}
    \caption{(Left) Labelling of the 4 quadrants. (Right) Rotations involved in step $A^5C^4BE$.}
    \label{fig:sample}
\end{figure}

Label the four quadrants (half-courts) $A, B, C, D$ as shown in Figure \ref{fig:sample}. 
The change in configuration, from the start of a game to the end of that game, can be represented symbolically in the form 
$$A^{x_1}C^{x_2}B^{x_3}D^{x_4},$$
where $0 \le x_i \le 5$ indicates the number of positions (modulo $6$) rotated by the team in the relevant quadrant. 
Figure  \ref{fig:notation} gives an illustration.

\setlength{\unitlength}{0.15in}
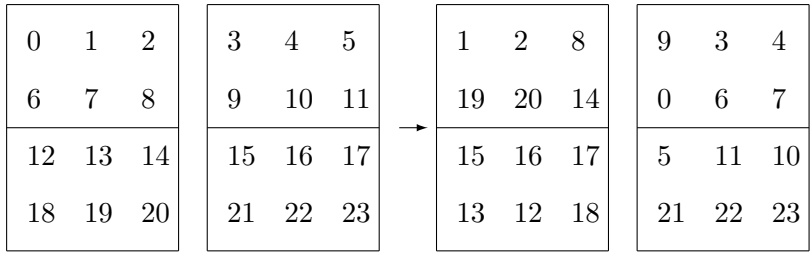
\begin{figure}[h]
\begin{picture}(0,10)(-2,0)
\put(0.3,0){\line(1,0){6}}
\put(0.3,4.3){\line(1,0){6}}
\put(0.3,8.6){\line(1,0){6}}
\put(0.3,0){\line(0,1){8.6}}
\put(6.3,0){\line(0,1){8.6}}
\put(1,1){18}
\put(3,1){19}
\put(5,1){20}
\put(1,3){12}
\put(3,3){13}
\put(5,3){14}
\put(1,5){6}
\put(3,5){7}
\put(5,5){8}
\put(1,7){0}
\put(3,7){1}
\put(5,7){2}
\put(7.3,0){\line(1,0){6}}
\put(7.3,4.3){\line(1,0){6}}
\put(7.3,8.6){\line(1,0){6}}
\put(7.3,0){\line(0,1){8.6}}
\put(13.3,0){\line(0,1){8.6}}
\put(8,1){21}
\put(10,1){22}
\put(12,1){23}
\put(8,3){15}
\put(10,3){16}
\put(12,3){17}
\put(8,5){9}
\put(10,5){10}
\put(12,5){11}
\put(8,7){3}
\put(10,7){4}
\put(12,7){5}
\put(15.3,0){\line(1,0){6}}
\put(15.3,4.3){\line(1,0){6}}
\put(15.3,8.6){\line(1,0){6}}
\put(15.3,0){\line(0,1){8.6}}
\put(21.3,0){\line(0,1){8.6}}
\put(16,1){13}
\put(18,1){12}
\put(20,1){18}
\put(16,3){15}
\put(18,3){16}
\put(20,3){17}
\put(16,5){19}
\put(18,5){20}
\put(20,5){14}
\put(16,7){1}
\put(18,7){2}
\put(20,7){8}
\put(22.3,0){\line(1,0){6}}
\put(22.3,4.3){\line(1,0){6}}
\put(22.3,8.6){\line(1,0){6}}
\put(22.3,0){\line(0,1){8.6}}
\put(28.3,0){\line(0,1){8.6}}
\put(23,1){21}
\put(25,1){22}
\put(27,1){23}
\put(23,3){5}
\put(25,3){11}
\put(27,3){10}
\put(23,5){0}
\put(25,5){6}
\put(27,5){7}
\put(23,7){9}
\put(25,7){3}
\put(27,7){4}
\put(14,4.3){\vector(1,0){1}}
\end{picture}
\caption{The effect of step $A^5C^4BE$.} 
\label{fig:notation}
\end{figure}

So the allowable values are 
$x_1 - x_2 \in \{-1,0,1\}$ and $x_3 - x_4 \in \{-1,0,1\}$ modulo $6$.
We then append a symbol $E$ to indicate the final movement 
(the front row players in each quadrant move to the same positions in the next quadrant).
This provides a coding of a step of the chain.  
For brevity we omit  any $x_i = 0$ term and write $B$ instead of $B^1$.
So a  typical step is coded in a format like $A^5C^4BE$.
The reader may check that the Figure \ref{Fig:1} example is $A^5C^4B^3D^2E$.

A sequence of steps can then be specified by concatenation: so  $A^5C^4BEEDE$ represents 3 steps of the chain, the second ($EE$) step indicating a game with zero (modulo 6) rotations of each team before the front row switch. 
We will name certain sequences later as $X, F, G, H$ in describing the construction.
The number of steps in a sequence is just the number of $E$'s, when expanded fully.
In writing the sequences (such as the definition of $X$ below) we often include spaces for visual clarity
but the spaces have no mathematical significance.

One aspect of this notation may be confusing.
The sequence $EEEE$ would code the identity move.  
That means that $BE\  EEEE$ has the same effect as $BE$.
But note that $BEEEE$ is different, and in fact will be a useful device because
it has the effect of rotating the players in quadrant $B$ while fixing all other players -- 
see Figure \ref{Fig:BE}. 
Note also that $B^5EEEE$ is the analogous back-rotation.
This syntax issue explains why we sometimes (e.g. in the definition of $F$ below) 
need to include an initial $EEEE$ in the definition.

\setlength{\unitlength}{0.15in}
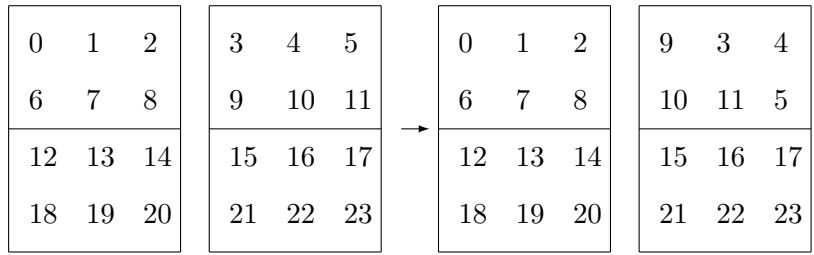
\begin{figure}[h]
\begin{picture}(0,10)(-2,0)
\put(0.3,0){\line(1,0){6}}
\put(0.3,4.3){\line(1,0){6}}
\put(0.3,8.6){\line(1,0){6}}
\put(0.3,0){\line(0,1){8.6}}
\put(6.3,0){\line(0,1){8.6}}
\put(1,1){18}
\put(3,1){19}
\put(5,1){20}
\put(1,3){12}
\put(3,3){13}
\put(5,3){14}
\put(1,5){6}
\put(3,5){7}
\put(5,5){8}
\put(1,7){0}
\put(3,7){1}
\put(5,7){2}
\put(7.3,0){\line(1,0){6}}
\put(7.3,4.3){\line(1,0){6}}
\put(7.3,8.6){\line(1,0){6}}
\put(7.3,0){\line(0,1){8.6}}
\put(13.3,0){\line(0,1){8.6}}
\put(8,1){21}
\put(10,1){22}
\put(12,1){23}
\put(8,3){15}
\put(10,3){16}
\put(12,3){17}
\put(8,5){9}
\put(10,5){10}
\put(12,5){11}
\put(8,7){3}
\put(10,7){4}
\put(12,7){5}
\put(15.3,0){\line(1,0){6}}
\put(15.3,4.3){\line(1,0){6}}
\put(15.3,8.6){\line(1,0){6}}
\put(15.3,0){\line(0,1){8.6}}
\put(21.3,0){\line(0,1){8.6}}
\put(16,1){18}
\put(18,1){19}
\put(20,1){20}
\put(16,3){12}
\put(18,3){13}
\put(20,3){14}
\put(16,5){6}
\put(18,5){7}
\put(20,5){8}
\put(16,7){0}
\put(18,7){1}
\put(20,7){2}
\put(22.3,0){\line(1,0){6}}
\put(22.3,4.3){\line(1,0){6}}
\put(22.3,8.6){\line(1,0){6}}
\put(22.3,0){\line(0,1){8.6}}
\put(28.3,0){\line(0,1){8.6}}
\put(23,1){21}
\put(25,1){22}
\put(27,1){23}
\put(23,3){15}
\put(25,3){16}
\put(27,3){17}
\put(23,5){10}
\put(25,5){11}
\put(27,5){5}
\put(23,7){9}
\put(25,7){3}
\put(27,7){4}
\put(14,4.3){\vector(1,0){1}}
\end{picture}
\caption{The effect of sequance $BEEEE$.} 
\label{Fig:BE}
\end{figure}

\subsection{High level description 1}
It is an elementary fact that any permutation of a card deck can be obtained by a sequence of transpositions of two adjacent cards. 
Indeed the ``random adjacent transposition" shuffling scheme is one of the original and most deeply studied examples in  the modern theory of mixing times 
\cite{me-RW,lacoin,wilson}.
By analogy, we start by showing that any transposition of two players on the same quadrant can be obtained by some sequence of steps.
There are 3 cases, depending on the initial distance between the two players,
illustrated in Figure \ref{fig:transpositions}, and we will exhibit sequences $F, G, H$ for each case.
We will show the first case (adjacent players) in detail.

\setlength{\unitlength}{0.15in}
\begin{figure}[h]
\begin{subfigure}[b]{\textwidth}
\begin{picture}(0,10)
\put(0.3,0){\line(1,0){6}}
\put(0.3,4.3){\line(1,0){6}}
\put(0.3,8.6){\line(1,0){6}}
\put(0.3,0){\line(0,1){8.6}}
\put(6.3,0){\line(0,1){8.6}}
\put(1,1){18}
\put(3,1){19}
\put(5,1){20}
\put(1,3){12}
\put(3,3){13}
\put(5,3){14}
\put(1,5){6}
\put(3,5){7}
\put(5,5){8}
\put(1,7){0}
\put(3,7){1}
\put(5,7){2}
\put(7.3,0){\line(1,0){6}}
\put(7.3,4.3){\line(1,0){6}}
\put(7.3,8.6){\line(1,0){6}}
\put(7.3,0){\line(0,1){8.6}}
\put(13.3,0){\line(0,1){8.6}}
\put(8,1){21}
\put(10,1){22}
\put(12,1){23}
\put(8,3){15}
\put(10,3){16}
\put(12,3){17}
\put(8,5){9}
\put(10,5){10}
\put(12,5){11}
\put(8,7){3}
\put(10,7){4}
\put(12,7){5}
\put(15.3,0){\line(1,0){6}}
\put(15.3,4.3){\line(1,0){6}}
\put(15.3,8.6){\line(1,0){6}}
\put(15.3,0){\line(0,1){8.6}}
\put(21.3,0){\line(0,1){8.6}}
\put(16,1){18}
\put(18,1){19}
\put(20,1){20}
\put(16,3){12}
\put(18,3){13}
\put(20,3){14}
\put(16,5){6}
\put(18,5){7}
\put(20,5){8}
\put(16,7){1}
\put(18,7){0}
\put(20,7){2}
\put(22.3,0){\line(1,0){6}}
\put(22.3,4.3){\line(1,0){6}}
\put(22.3,8.6){\line(1,0){6}}
\put(22.3,0){\line(0,1){8.6}}
\put(28.3,0){\line(0,1){8.6}}
\put(23,1){21}
\put(25,1){22}
\put(27,1){23}
\put(23,3){15}
\put(25,3){16}
\put(27,3){17}
\put(23,5){9}
\put(25,5){10}
\put(27,5){11}
\put(23,7){3}
\put(25,7){4}
\put(27,7){5}
\put(14,4.3){\vector(1,0){1}}
\end{picture}
\caption{Transpose two adjacent players in the same quadrant (sequence F)} \label{fig:transpose2}
\end{subfigure}

\begin{subfigure}[b]{\textwidth}
\begin{picture}(0,10)
\put(0.3,0){\line(1,0){6}}
\put(0.3,4.3){\line(1,0){6}}
\put(0.3,8.6){\line(1,0){6}}
\put(0.3,0){\line(0,1){8.6}}
\put(6.3,0){\line(0,1){8.6}}
\put(1,1){18}
\put(3,1){19}
\put(5,1){20}
\put(1,3){12}
\put(3,3){13}
\put(5,3){14}
\put(1,5){6}
\put(3,5){7}
\put(5,5){8}
\put(1,7){0}
\put(3,7){1}
\put(5,7){2}
\put(7.3,0){\line(1,0){6}}
\put(7.3,4.3){\line(1,0){6}}
\put(7.3,8.6){\line(1,0){6}}
\put(7.3,0){\line(0,1){8.6}}
\put(13.3,0){\line(0,1){8.6}}
\put(8,1){21}
\put(10,1){22}
\put(12,1){23}
\put(8,3){15}
\put(10,3){16}
\put(12,3){17}
\put(8,5){9}
\put(10,5){10}
\put(12,5){11}
\put(8,7){3}
\put(10,7){4}
\put(12,7){5}
\put(15.3,0){\line(1,0){6}}
\put(15.3,4.3){\line(1,0){6}}
\put(15.3,8.6){\line(1,0){6}}
\put(15.3,0){\line(0,1){8.6}}
\put(21.3,0){\line(0,1){8.6}}
\put(16,1){18}
\put(18,1){19}
\put(20,1){20}
\put(16,3){12}
\put(18,3){13}
\put(20,3){14}
\put(16,5){6}
\put(18,5){7}
\put(20,5){8}
\put(16,7){2}
\put(18,7){1}
\put(20,7){0}
\put(22.3,0){\line(1,0){6}}
\put(22.3,4.3){\line(1,0){6}}
\put(22.3,8.6){\line(1,0){6}}
\put(22.3,0){\line(0,1){8.6}}
\put(28.3,0){\line(0,1){8.6}}
\put(23,1){21}
\put(25,1){22}
\put(27,1){23}
\put(23,3){15}
\put(25,3){16}
\put(27,3){17}
\put(23,5){9}
\put(25,5){10}
\put(27,5){11}
\put(23,7){3}
\put(25,7){4}
\put(27,7){5}
\put(14,4.3){\vector(1,0){1}}
\end{picture}
\caption{Transpose two players in the same quadrant with one space in between (G)}
\label{fig:transpose3}
\end{subfigure}

\begin{subfigure}[b]{\textwidth}
\begin{picture}(0,10)
\put(0.3,0){\line(1,0){6}}
\put(0.3,4.3){\line(1,0){6}}
\put(0.3,8.6){\line(1,0){6}}
\put(0.3,0){\line(0,1){8.6}}
\put(6.3,0){\line(0,1){8.6}}
\put(1,1){18}
\put(3,1){19}
\put(5,1){20}
\put(1,3){12}
\put(3,3){13}
\put(5,3){14}
\put(1,5){6}
\put(3,5){7}
\put(5,5){8}
\put(1,7){0}
\put(3,7){1}
\put(5,7){2}
\put(7.3,0){\line(1,0){6}}
\put(7.3,4.3){\line(1,0){6}}
\put(7.3,8.6){\line(1,0){6}}
\put(7.3,0){\line(0,1){8.6}}
\put(13.3,0){\line(0,1){8.6}}
\put(8,1){21}
\put(10,1){22}
\put(12,1){23}
\put(8,3){15}
\put(10,3){16}
\put(12,3){17}
\put(8,5){9}
\put(10,5){10}
\put(12,5){11}
\put(8,7){3}
\put(10,7){4}
\put(12,7){5}
\put(15.3,0){\line(1,0){6}}
\put(15.3,4.3){\line(1,0){6}}
\put(15.3,8.6){\line(1,0){6}}
\put(15.3,0){\line(0,1){8.6}}
\put(21.3,0){\line(0,1){8.6}}
\put(16,1){18}
\put(18,1){19}
\put(20,1){20}
\put(16,3){12}
\put(18,3){13}
\put(20,3){14}
\put(16,5){6}
\put(18,5){7}
\put(20,5){0}
\put(16,7){8}
\put(18,7){1}
\put(20,7){2}
\put(22.3,0){\line(1,0){6}}
\put(22.3,4.3){\line(1,0){6}}
\put(22.3,8.6){\line(1,0){6}}
\put(22.3,0){\line(0,1){8.6}}
\put(28.3,0){\line(0,1){8.6}}
\put(23,1){21}
\put(25,1){22}
\put(27,1){23}
\put(23,3){15}
\put(25,3){16}
\put(27,3){17}
\put(23,5){9}
\put(25,5){10}
\put(27,5){11}
\put(23,7){3}
\put(25,7){4}
\put(27,7){5}
\put(14,4.3){\vector(1,0){1}}
\end{picture}
\caption{Transpose two players in the same quadrant with two spaces in between (H)}
\label{fig:transpose4}
\end{subfigure}
\caption{Transpositions achieved by the specific sequences $F, G, H$.}\label{fig:transpositions}
\end{figure}

\newpage
\subsection{Sequences that transpose two players}
\label{sec:transpose}
We start by  introducing a 16 step sequence  $X$ defined as
$$(X:=) \quad AE\;\;B^2D^3EEE\;\;A^2C^3E\;\;B^3D^3EEE\;\;
A^5E\;\;B^5EEE\;\;AE\;\;B^5EEE\;\;A C^3.$$
The step-by-step trajectory of sequence $X$ is shown in Figure \ref{fig:YYY} below, 
which demonstrates that the effect of $X$ is as shown in Figure \ref{fig:X}. 
The introduction of this $X$ is somewhat magical and hard to explain, but note that for some players it is like a reverse step of the chain.

\setlength{\unitlength}{0.15in}
\begin{figure}[h]
\begin{picture}(0,10)
\put(0.3,0){\line(1,0){6}}
\put(0.3,4.3){\line(1,0){6}}
\put(0.3,8.6){\line(1,0){6}}
\put(0.3,0){\line(0,1){8.6}}
\put(6.3,0){\line(0,1){8.6}}
\put(1,1){18}
\put(3,1){19}
\put(5,1){20}
\put(1,3){12}
\put(3,3){13}
\put(5,3){14}
\put(1,5){6}
\put(3,5){7}
\put(5,5){8}
\put(1,7){0}
\put(3,7){1}
\put(5,7){2}
\put(7.3,0){\line(1,0){6}}
\put(7.3,4.3){\line(1,0){6}}
\put(7.3,8.6){\line(1,0){6}}
\put(7.3,0){\line(0,1){8.6}}
\put(13.3,0){\line(0,1){8.6}}
\put(8,1){21}
\put(10,1){22}
\put(12,1){23}
\put(8,3){15}
\put(10,3){16}
\put(12,3){17}
\put(8,5){9}
\put(10,5){10}
\put(12,5){11}
\put(8,7){3}
\put(10,7){4}
\put(12,7){5}
\put(15.3,0){\line(1,0){6}}
\put(15.3,4.3){\line(1,0){6}}
\put(15.3,8.6){\line(1,0){6}}
\put(15.3,0){\line(0,1){8.6}}
\put(21.3,0){\line(0,1){8.6}}
\put(16,1){18}
\put(18,1){19}
\put(20,1){20}
\put(16,3){12}
\put(18,3){13}
\put(20,3){14}
\put(16,5){5}
\put(18,5){3}
\put(20,5){6}
\put(16,7){1}
\put(18,7){0}
\put(20,7){2}
\put(22.3,0){\line(1,0){6}}
\put(22.3,4.3){\line(1,0){6}}
\put(22.3,8.6){\line(1,0){6}}
\put(22.3,0){\line(0,1){8.6}}
\put(28.3,0){\line(0,1){8.6}}
\put(23,1){21}
\put(25,1){22}
\put(27,1){23}
\put(23,3){15}
\put(25,3){16}
\put(27,3){17}
\put(23,5){9}
\put(25,5){10}
\put(27,5){11}
\put(23,7){4}
\put(25,7){7}
\put(27,7){8}
\put(14,4.3){\vector(1,0){1}}
\end{picture}
\caption{The effect of sequence $X$.} 
\label{fig:X}
\end{figure}
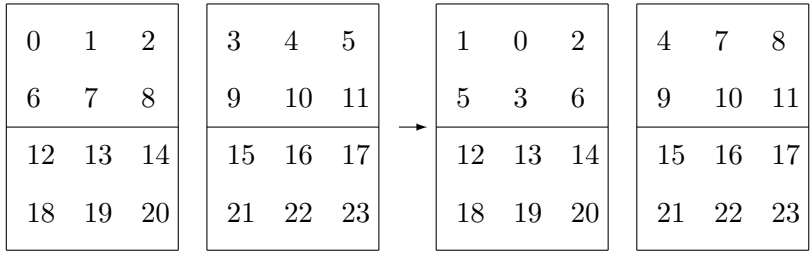

\vspace*{-1.3in}

\newcommand{\move}[2]{
\setsepchar{ }
\readlist\arg{#1}
\begin{subfigure}[b]{5.6cm}
\begin{tikzpicture}[scale=0.70]

\draw (-0.2,-1.5) -- (2.2,-1.5);
\draw (-0.2,0.5) -- (-0.2,-3.5);
\draw (-0.2,0.5) -- (2.2,0.5);
\draw (-0.2,-3.5) -- (2.2,-3.5);
\draw (2.2,0.5) -- (2.2,-3.5); 
 
\draw (3,-1.5) -- (5.4,-1.5);
\draw (3,0.5) -- (3,-3.5);
\draw (3,0.5) -- (5.4,0.5);
\draw (3,-3.5) -- (5.4,-3.5);
\draw (5.4,0.5) -- (5.4,-3.5); 

\draw   (0,1.1) node[text width=6cm] {#2};

\node at (0.2,-0.1) {\arg[1]} ; 
\node at (1.0,-0.1) {\arg[2]} ; 
\node at (1.8,-0.1) {\arg[3]} ; 
\node at (3.4,-0.1) {\arg[4]} ; 
\node at (4.2,-0.1) {\arg[5]} ; 
\node at (5.0,-0.1) {\arg[6]} ; 
\node at (0.2,-0.9) {\arg[7]} ; 
\node at (1.0,-0.9) {\arg[8]} ; 
\node at (1.8,-0.9) {\arg[9]} ; 
\node at (3.4,-0.9) {\arg[10]} ; 
\node at (4.2,-0.9) {\arg[11]} ; 
\node at (5.0,-0.9) {\arg[12]} ; 

\node at (0.2,-2.1) {\arg[13]} ; 
\node at (1.0,-2.1) {\arg[14]} ; 
\node at (1.8,-2.1) {\arg[15]} ; 
\node at (3.4,-2.1) {\arg[16]} ; 
\node at (4.2,-2.1) {\arg[17]} ; 
\node at (5.0,-2.1) {\arg[18]} ; 
\node at (0.2,-2.9) {\arg[19]} ; 
\node at (1.0,-2.9) {\arg[20]} ; 
\node at (1.8,-2.9) {\arg[21]} ; 
\node at (3.4,-2.9) {\arg[22]} ; 
\node at (4.2,-2.9) {\arg[23]} ; 
\node at (5.0,-2.9) {\arg[24]} ; 

 \end{tikzpicture}
 \end{subfigure}
 \hspace{0.5cm}
}

\begin{figure}
    \centering
    \move{6 0 1 3 4 5 14 13 12 7 8 2 15 16 17 11 10 9 18 19 20 21 22 23}{1. $AE$ (Separate player $2$ from $0$ and $1$.)}
    \move{6 0 1 8 7 3 2 5 4 21 22 23 12 13 14 15 16 17 18 19 20 9 10 11}{2 to 4. $B^2D^3EEE$ (Move player $2$ to the front right position in A.)}
    \move{5 2 6 8 7 3 18 19 20 4 1 0 15 16 17 23 22 21 14 13 12 9 10 11}{5. $A^2C^3E$ (Separate players $0$ and $1$ from $2$.)}
    \move{5 2 6 0 1 4 3 7 8 9 10 11 20 19 18 15 16 17 14 13 12 21 22 23}{6 to 8. $B^3D^3EEE$ (Move player $0$ at the back right of B.)}
    \move{2 6 8 0 1 4 18 19 20 5 3 7 15 16 17 11 10 9 14 13 12 21 22 23}{9. $A^5E$ (Rotate  $A$ to move player $2$ at the back right of A.)}
    \move{2 6 8 1 4 7 0 5 3 9 10 11 20 19 18 15 16 17 14 13 12 21 22 23}{10 to 12. $B^5EEE$ (Move player $2$ to the left of $0$.) }
    \move{0 2 6 1 4 7 18 19 20 5 3 8 15 16 17 11 10 9 14 13 12 21 22 23}{13. $AE$}
    \move{0 2 6 4 7 8 1 5 3 9 10 11 20 19 18 15 16 17 14 13 12 21 22 23}{14 to 16. $B^5EEE$ (Move player 1 to the left of $02$.)}
    \move{1 0 2 4 7 8 5 3 6 9 10 11 12 13 14 15 16 17 18 19 20 21 22 23}{17 to 20. $AC^2EEEEC$ (Same effect as $AC^3$)}

    \caption{Step-by-step trajectory of sequence $X$.}
    \label{fig:YYY}
\end{figure}

\newpage
We can now define the sequence $F$ that transposes the two adjacent players at the left corner of the back row of quadrant $A$, as shown in Figure \ref{fig:transpose2}.
Essentially it is just 3 applications of $X$.
Precisely
\begin{align*}
(F:=) \quad &EEEEX\;\;X\;\;XEEEE.
\end{align*}
So $F$ involves 
$68$ steps of the chain.
Figure \ref{fig:ZZZ} shows the step-by-step trajectory of sequence $F$.

\vspace{0.3in}

\begin{figure}
    \centering
    \move{1 0 2 4 7 8 5 3 6 9 10 11 12 13 14 15 16 17 18 19 20 21 22 23}{1. $EEEEX$}
    \move{0 1 2 7 3 6 8 4 5 9 10 11 12 13 14 15 16 17 18 19 20 21 22 23}{2. $X$}
    \move{1 0 2 3 4 5 6 7 8 9 10 11 12 13 14 15 16 17 18 19 20 21 22 23}{3. $XEEEE$}
    \caption{Step-by-step trajectory of sequence $F$.}
    \label{fig:ZZZ}
\end{figure}

\paragraph{The other transposition sequences.}
To transpose the two players at the back row of quadrant $A$ with one space in between, as shown in Figure \ref{fig:transpose3}, we use the sequence $G$ defined as
$$(G:=) \quad EEEE\;\;A^5FA\;\;F\;\;A^5FA\;\;EEEE.$$
This works because the effect of these sequences is to alter the back row as
\[ 012 \to 021 \to 201 \to 210 .\]
Finally, to transpose the players at the upper left corner and at the lower right corner in quadrant $A$, as shown in Figure \ref{fig:transpose4}, 
we use the sequence $H$ defined as 
$$(H:=) \quad FA^5\;\;FA^5\;\;FAFAF.$$
The reader may check that this works, and is one place where the initial $EEEE$ in the definition of $F$ is needed.

\subsection{High level description 2}
By symmetry, to prove irreducibility it is enough to prove that, from any initial state, one can reach 
(by some sequence of allowable steps)
 the reference state shown in Figure \ref{fig:goal}, where each player  $i$ is in position $i$.
 Because any permutation of the 6 players in quadrant $A$ can be derived from a sequence of transpositions,
 from the existence of transposition sequences (section \ref{sec:transpose}) it suffices to show that we can move,
 by some sequence of steps, player $i$ to position $i$, for every player in all the other quadrants $B, C, D$.
 We will move players to positions row by row, in the following order
 \begin{itemize}
 \item back row of $D$: (21, 22, 23)
 \item back row of $C$: (18, 19, 20)
 \item front row of $D$: (17, 16, 15) 
 \item front row of $C$: (14, 13, 12) 
 \item back row of $B$: (5, 4, 3) 
 \item front row of $B$: (9, 10, 11)
 \end{itemize}
 Each row in turn is {\em fixed}, in that it remains in place after each subsequent row has been moved to its position.

\setlength{\unitlength}{0.15in}
\begin{figure}[h]
\centering
\begin{picture}(18,10)
\put(0.3,0){\line(1,0){6}}
\put(0.3,4.3){\line(1,0){6}}
\put(0.3,8.6){\line(1,0){6}}
\put(0.3,0){\line(0,1){8.6}}
\put(6.3,0){\line(0,1){8.6}}
\put(1,1){18}
\put(3,1){19}
\put(5,1){20}
\put(1,3){12}
\put(3,3){13}
\put(5,3){14}
\put(1,5){6}
\put(3,5){7}
\put(5,5){8}
\put(1,7){0}
\put(3,7){1}
\put(5,7){2}
\put(13,1){21}
\put(15,1){22}
\put(17,1){23}
\put(13,3){15}
\put(15,3){16}
\put(17,3){17}
\put(13,5){9}
\put(15,5){10}
\put(17,5){11}
\put(13,7){3}
\put(15,7){4}
\put(17,7){5}
\put(12.3,0){\line(1,0){6}}
\put(12.3,4.3){\line(1,0){6}}
\put(12.3,8.6){\line(1,0){6}}
\put(12.3,0){\line(0,1){8.6}}
\put(18.3,0){\line(0,1){8.6}}
\put(2,6){A}
\put(14,6){B}
\put(2,2){C}
\put(14,2){D}
\end{picture}
\caption{Reference State}  
\label{fig:goal}
\end{figure}
In the next section we will describe the algorithm in words.
 A key point is that we will use quadrant $A$ as a kind of temporary stopover for players in transit.

 \subsection{The algorithm}
 \label{sec:algorithm}
 
 As noted in Figure  \ref{Fig:BE} a sequence like $BEEE$ has the effect of rotating a given quadrant by one position,
 and one can repeat such a sequence.  So we can use phrases such as 
 ``{\em rotate} player $i$ to position $j$" (in the same quadrant). 
 We will call $E$ the {\em migration} step, so by repeating step $E$ we can use ``migrate player $i$ to position $j$" (where player $i$ has the same relative position (e.g. front, right) as $j$, in the front row.  
 (But remember that all 12 front-row players migrate.)
 In both cases we of course only do the move if necessary, that is if the player is not already in the desired position.

 The algorithm is based on variations of the following ``Procedure P", where $P$ is one of the quadrants 
 $B, C, D$, and where the procedure  acts to move the required 3 players to the back row of $P$.
 
 \vspace{0.1in}  \noindent 
{\bf Procedure P.}

\vspace{0.1in}  \noindent  ({\bf 1}) 
Label the back row positions of $P$ as $a, b, c$, left-to-right.

(So for $P = D$ we have $(a,b,c) = (21, 22, 23)$).

\vspace{0.1in}  \noindent  ({\bf 2}) 
Rotate (if in back row) player $a$ to front row. Migrate player $a$ to quadrant $A$ and rotate to front right position.
Migrate player $a$ to the front right position  in  $P$.

\vspace{0.1in}  \noindent  ({\bf 3}) 
Rotate player $a$ to back right position  in  $P$.

\vspace{0.1in}  \noindent  ({\bf 4a}) 
 If player $b$ is not in  $P$, repeat actions ({\bf 2}) for player $b$. This moves player $b$ to 
 the front right position  in  $P$.
 
\vspace{0.1in}  \noindent  ({\bf 4b}) 
 Else: rotate  $P$ so that player $b$ is in front row while player $a$ is in back row.
Migrate player $b$ to quadrant $A$.  Rotate $A$ to move player $b$ to front right.
Rotate $P$ so that player $a$ is returned to back right position  in  $P$.
Migrate player $b$ to the front right position  in  $P$.

\vspace{0.1in}  \noindent  ({\bf 5}) 
Now players $(a,b)$ are in (back right, front right) positions in $P$: rotate one space to
(back center, back right) positions.

(Next we move player $c$ in essentially the same way as player $b$,  Here are the details.)

\vspace{0.1in}  \noindent  ({\bf 6a}) 
 If player $c$ is not in  $P$, repeat actions ({\bf 2}) for player $c$. This moves player $c$ to 
 the front right position  in  $P$.
 
\vspace{0.1in}  \noindent  ({\bf 6b}) 
 Else: rotate  $P$ so that player $c$ is in front row while players $a,b$ are in back row.
Migrate player $c$ to quadrant $A$.  Rotate $A$ to move player $c$ to front right.
Rotate $P$ so that players $(a,b)$ are returned to (back center, back right) positions in  $P$.
Migrate player $c$ to the front right position  in  $P$.

\vspace{0.1in}  \noindent  ({\bf 7}) 
Rotate players $(a,b,c)$ to (back left, back center, back right) positions in $P$.

\vspace{0.1in}  \noindent  
 {\bf end Procedure P.}
 
 \medskip \noindent
 Procedure $P$ moves the required 3 players to the back row of $P$.
 In the algorithm below, 
 once a back row is placed, it is ``fixed" and those players are never moved subsequently. 
 This is because that would require a rotation of $P$ at some subsequent use of  ({\bf 2}), which cannot happen 
 because of the   {\em if in back row} condition in ({\bf 2}): these players are already in place.
 
\bigskip \noindent
{\bf The algorithm.}

\vspace{0.1in}  \noindent  ({\bf 8}) 
Apply Procedure $D$, then Procedure $C$.

(This fixes the back rows of $C$ and $D$.  
For the third row (the front row in $D$) we will use a little trick, to first arrange them on a back row.)

\vspace{0.1in}  \noindent  ({\bf 9}) 
Here, we consider players (17,16,15)  who need to be moved to the front row of $D$.
But we label them as $(a,b,c)$ = (5,4,3)  and use Procedure $B$ to move them to the back row of $B$. Then revert to labels $(17,16,15)$.

\vspace{0.1in}  \noindent  ({\bf 10}) 
Rotate players (17,16,15) three positions so they become the front row of $B$.

\vspace{0.1in}  \noindent  ({\bf 11}) 
Migrate one step, so  players (17,16,15)  become the front row in $D$.

(The previous back rows remain fixed. Next we consider the fourth row, the front row of $C$. 
Here another small complication arises; any migration will move the front row of $D$  (fixed above), so we need to ensure that it will not be rotated (we can migrate if necessary), and it is migrated back into place before finishing.)

\vspace{0.1in}  \noindent  ({\bf 12}) 
Here we are considering players (14,13,12).
As in ({\bf 9}), label them as $(a,b,c) = (5, 4, 3)$  and use Procedure $B$ to move them to the back row of $B$.
Then revert to labels (14,13,12).
 This involves a certain number of migrates. When applying Procedure $B$, we make sure players $(17, 16, 15)$ remain at the front row in some quadrant at all times by not rotating them, and the fixed back court players (18-23) are not moved.
 
 (If we need to rotate a quadrant where players (17,16,15) are in the front row, migrate first until that front row does not contain players $(17, 16, 15)$, $b$ if $b$ is not at the back row of $P$, and $c$ if $b$ is at the back row of $P$ and $c$ is not, then rotate that quadrant.)

\vspace{0.1in}  \noindent  ({\bf 13}) 
Migrate players (17, 16, 15) to make them the front row in quadrant $A$. 
Rotate players (14,13,12) in $B$ three positions so they become the front row in $B$.
Migrate two turns.

(Now all players 12-23 in quadrants $C$ and $D$ have been moved to their reference positions.)

\vspace{0.1in}  \noindent  ({\bf 14}) 
Apply Procedure $B$.  This moves players (5,4,3) to the back row of $B$, as required.
It involves some number of migrates, but the players (12-17) assigned to front rows of $C$ and $D$ 
should remain as front rows of adjacent quadrants, so make sure not to rotate them and migrate them back to their required positions.

(If a rotation is needed, migrate first until that front row does not contain players (12-17), $b$ if $b$ is not at the back row of $P$, and $c$ if $b$ is at the back row of $P$ and $c$ is not.)

(Now we have fixed all rows of $B, C, D$ except the front row of $B$. 
For the front row of $B$ we do a quite different scheme, using the 
transposition sequences from section \ref{sec:transpose}.)

\vspace{0.1in}  \noindent  ({\bf 15}) 
Here we are considering players (11,10,9).
Our first goal is to move them to target positions (2,1,0) in the  back row of $A$.
Any of those players (11,10,9) in $A$ can be moved to their target position via transpositions. 
Remaining players amongst (11,10,9) must be in the front row of $B$.
Migrate the front row of $B$ to the front row of $A$, transpose relevant players to target positions, and migrate back.

\vspace{0.1in}  \noindent  ({\bf 16}) 
 Players (11,10,9)  are at positions (2,1,0) in the  back row of $A$.  Migrate (17, 16, 15) to the front row of B. Rotate players (11,10,9) three positions so they become the front row of $A$, and migrate to front row of $B$.
 
 (Now all rows of $B,C,D$ are fixed.)

 \vspace{0.1in}  \noindent  ({\bf 17}) 
As noted earlier, any permutation of the 6 players in quadrant $A$ can be derived from a sequence of transpositions, so
 we can fix quadrant $A$ by using the transposition sequences $F, G, H$ from section \ref{sec:transpose}.
 
 \vspace{0.1in}  \noindent 
 {\bf End Algorithm}

\bigskip
This completes our proof of irreducibility.  Recall that the basic idea was to move players to positions row-by-row.  It seems likely that there is some alternate 
``basic idea" that would lead to a simpler proof, and so we have not tried to optimize the implementation of our basic idea.

\section{The big chain is aperiodic}
\label{sec:aperiodic}

We know that it is possible to move from one state to that same state in four steps - $EEEE$. 
So to prove that the big chain is aperiodic, it is sufficient to exhibit a sequence of $n$ steps, for some odd $n$, that move the reference state 
(Figure \ref{fig:goal}) to that same state.
Such a sequence is shown in Figure \ref{fig:aperiodic},
There are 163 $E$'s in the non-expanded notation, and 68, 424, and 340 
steps from the 1 $F$, 2 $G$'s, and 1 $H$, respectively, so there are 
995 steps in this sequence.

\newcommand{\movex}[2]{
\setsepchar{ }
\readlist\arg{#1}
\begin{subfigure}[b]{5.5cm}
\begin{tikzpicture}[scale=0.70]

\draw (-0.2,-1.5) -- (2.2,-1.5);
\draw (-0.2,0.5) -- (-0.2,-3.5);
\draw (-0.2,0.5) -- (2.2,0.5);
\draw (-0.2,-3.5) -- (2.2,-3.5);
\draw (2.2,0.5) -- (2.2,-3.5); 
 
\draw (3,-1.5) -- (5.4,-1.5);
\draw (3,0.5) -- (3,-3.5);
\draw (3,0.5) -- (5.4,0.5);
\draw (3,-3.5) -- (5.4,-3.5);
\draw (5.4,0.5) -- (5.4,-3.5); 

\draw   (0,2) node[text width=5.6cm] {#2};

\node at (0.2,-0.1) {\arg[1]} ; 
\node at (1.0,-0.1) {\arg[2]} ; 
\node at (1.8,-0.1) {\arg[3]} ; 
\node at (3.4,-0.1) {\arg[4]} ; 
\node at (4.2,-0.1) {\arg[5]} ; 
\node at (5.0,-0.1) {\arg[6]} ; 
\node at (0.2,-0.9) {\arg[7]} ; 
\node at (1.0,-0.9) {\arg[8]} ; 
\node at (1.8,-0.9) {\arg[9]} ; 
\node at (3.4,-0.9) {\arg[10]} ; 
\node at (4.2,-0.9) {\arg[11]} ; 
\node at (5.0,-0.9) {\arg[12]} ; 

\node at (0.2,-2.1) {\arg[13]} ; 
\node at (1.0,-2.1) {\arg[14]} ; 
\node at (1.8,-2.1) {\arg[15]} ; 
\node at (3.4,-2.1) {\arg[16]} ; 
\node at (4.2,-2.1) {\arg[17]} ; 
\node at (5.0,-2.1) {\arg[18]} ; 
\node at (0.2,-2.9) {\arg[19]} ; 
\node at (1.0,-2.9) {\arg[20]} ; 
\node at (1.8,-2.9) {\arg[21]} ; 
\node at (3.4,-2.9) {\arg[22]} ; 
\node at (4.2,-2.9) {\arg[23]} ; 
\node at (5.0,-2.9) {\arg[24]} ; 

 \end{tikzpicture}
 \end{subfigure}
 \hspace{0.5cm}
}

\begin{figure}
    \centering
    \movex{14 7 6 3 4 5 0 21 22 8 2 1 23 15 16 18 11 10 19 20 9 13 12 17}{$1\text{ to }13. AEEEE\;AEAE\;CE$ $DEEEE\;DEDE$ (Move to a certain other arrangement)}
    \movex{14 7 6 3 4 5 8 2 1 10 11 18 17 12 13 0 20 19 9 16 15 21 22 23}{ $14\text{ to }32. EEDEEEE\;DEEEE$ $DEE\;CEEEE\;CEEE\;D$ (Fixing back row of D)}
    \movex{11 10 14 3 4 5 9 16 15 8 2 1 17 12 0 13 7 6 18 19 20 21 22 23}{$33\text{ to }57. CEEEE\;CEEEE\;CEEE$ $AEEEE\;AEEE\;CEEE\;CEEEE\;C$ (Fixing back row of C)}
    \movex{9 12 0 3 4 11 6 7 13 14 10 5 1 2 8 15 16 17 18 19 20 21 22 23}{$58\text{ to }67. EAEEEE\;AEBEE\;AE\;BBBBBE$ (Fixing (17,16,15) - using Procedure B then some migration)}
    \movex{0 7 6 1 2 8 4 9 3 10 5 11 12 13 14 15 16 17 18 19 20 21 22 23}{$68\text{ to }87. BAEEEE\;AEBEEE\;AEBE$ $BEEEE\;BEEEE\;BEE$ (Fixing (14,13,12) - using Procedure B then some migration)}
    \movex{9 0 7 3 4 5 1 6 10 11 8 2 12 13 14 15 16 17 18 19 20 21 22 23}{$88\text{ to }99.BEEEE\;BEBEEE\;AEBEEE$  (Fixing (5,4,3) using Procedure B)}
    \movex{0 2 8 3 4 5 6 1 7 9 10 11 12 13 14 15 16 17 18 19 20 21 22 23}{$100\text{ to }135.AEEEE\;AEEEE\;AGAEEEE$ $AEEEE\;AEEEE\;AEEE\;AEEEE$ $AEEEE\;AEEEE\;AE$ (Fixing (9,10,11) using Steps 15 and 16 of \textbf{The Algorithm}}
    \movex{0 1 2 3 4 5 6 7 8 9 10 11 12 13 14 15 16 17 18 19 20 21 22 23}{$136\text{ to }163. AEEEE\;AHA\;AGA\;AEEEE$ $AEEEE\;AEEEE\;AFAEEEE\;AEEEE$ $AEEEE$ (Fixing A)}
    \caption{Sample Moves}
    \label{fig:aperiodic}
\end{figure}

\section{The friend chain}
\label{sec:friend}
Perhaps the most natural observable to study concerns the positions of two players, say {\em ego} and {\em friend}.  
Naively this would require a $24 \times 23$ state chain, but we can exploit some symmetries to make a 26 state chain indicating {\em relative} positions of the two players.
Doing so requires some care; the states are indicated in Figure \ref{Fig:26state}, as explained next.

First note that our process is not invariant under a quarter-turn of the 4 quadrants\footnote{A friend on the same team in quadrant $B$ might be an opponent in the next game; this cannot happen on quadrant $D$.}, but is invariant under a half-turn, so as in Figure \ref{Fig:26state}, we can assume {\em ego} is in the left court.

\setlength{\unitlength}{0.20in}
\begin{figure}[h]
\begin{picture}(18,28)(0,-8)
\put(0.3,0){\line(1,0){6.7}}
\put(0.3,4.3){\line(1,0){6.7}}
\put(0.3,8.6){\line(1,0){6.7}}
\put(0.3,0){\line(0,1){8.6}}
\put(7,0){\line(0,1){8.6}}
\put(1,1){1T5}
\put(3,1){1T6}
\put(5,1){$\bullet$}
\put(1,3){1T4}
\put(3,3){1T3}
\put(5,3){1T2}
\put(1,5){1O2}
\put(3,5){1O3}
\put(5,5){1O4}
\put(1,7){1O1}
\put(3,7){1O6}
\put(5,7){1O5}
\put(18,5){{\em ego} $\bullet$ in first quadrant}
\put(9.3,0){\line(1,0){6.7}}
\put(9.3,4.3){\line(1,0){6.7}}
\put(9.3,8.6){\line(1,0){6.7}}
\put(9.3,0){\line(0,1){8.6}}
\put(16,0){\line(0,1){8.6}}
\put(11.5,2){all 1-}
\put(11.5,6){all 1+}
\put(0.3,10){\line(1,0){6.7}}
\put(0.3,14.3){\line(1,0){6.7}}
\put(0.3,18.6){\line(1,0){6.7}}
\put(0.3,10){\line(0,1){8.6}}
\put(7,10){\line(0,1){8.6}}
\put(1,11){2O3}
\put(3,11){2O4}
\put(5,11){2O5}
\put(1,13){2O2}
\put(3,13){2O1}
\put(5,13){2O6}
\put(1,15){2T6}
\put(3,15){$\bullet$}
\put(5,15){2T2}
\put(1,17){2T5}
\put(3,17){2T4}
\put(5,17){2T3}
\put(18,15){{\em ego}  $\bullet$  in second quadrant}
\put(9.3,10){\line(1,0){6.7}}
\put(9.3,14.3){\line(1,0){6.7}}
\put(9.3,18.6){\line(1,0){6.7}}
\put(9.3,10){\line(0,1){8.6}}
\put(16,10){\line(0,1){8.6}}

\put(11.5,12){all 2-}
\put(11.5,16){all 2+}
\end{picture}
\vspace*{-1.5in}
\caption{Positions of {\em friend} relative to {\em ego} $\bullet$.}  
\label{Fig:26state}
\end{figure}

The states of what we will call the {\em friend chain} indicate relative positions at the start of a game.  
To describe the state, first record which quadrant {\em ego} is in (denoted here as initial 1 or 2; that is, quadrants $C$ or $A$) and then whether {\em friend} is in the same team (denoted T) or the current opposing team (denoted O), or on the other court.
If on the other court, we only need to note which quadrant (because the position will be randomized  during the game), so denote by one of (1+, 1-, 2+, 2-) as illustrated.
Finally, writing temporarily ego* for the opponent in the same position as ego, we indicate a friend's position as 1, 2, 3, 4, 5, or 6, counter-clockwise from ego or ego*.

This adds up to 26 states, and  one can check that the ``big" chain rules define a Markov chain on these states, with transition matrix $\bP$ as shown in Figures \ref{Fig:3} and \ref{Fig:4}.
Its stationary distribution $\pi$  is induced from the uniform stationary distribution of the ``big" chain:
probability $\tfrac6{46}$ for each of the states $1+, 1-, 2+, 2-$ and probability $\tfrac1{46}$ for each of the 22 remaining states.

\begin{figure}
$
\begin{array}{c|ccccccccccccc}
&1+ & 1- & 1T2 & 1T3 & 1T4 & 1T5 & 1T6 & 1O1 & 1O2 & 1O3 & 1O4 & 1O5 & 1O6 \\
\hline
1+ & 1/4 & 1/4 &&&&&&&&&&\\
1- & & 1/4 & 1/36 & 2/36 & 3/36 & 2/36 & 1/36 &&&&&&\\
1T2 & && 1/3 &&&&&& 1/6 &&&&\\
1T3 & &&& 1/6 &&&&&& 1/3 &&&\\
1T4 & &&&  &&&&&&& 1/2 &&\\
1T5 &&&&& & 1/6 &&&&&& 1/3 &\\
1T6 &&&&& && 1/3 &&&&&& 1/6 \\
1O1 &1/12 &&  &&&&&1/4& 1/12 &&&&1/12\\
1O2 &1/6 &&  &&&&&3/24& 1/6&1/24&&&\\
1O3 &1/3 &&  &&&&&& 1/12 &1/12&&&\\
1O4 &5/12 &&  &&&&&&&1/24& &1/24&\\
1O5&1/3 &&  &&&&&&&&& 1/12 &1/12\\
1O6&1/6 &&  &&&&&3/24&&&& 1/24&1/6\\

2T2 &&1/6&1/3&&&&&&&& &  & \\
2T3 &&1/3&&1/6&&&&&&& &  & \\
2T4 &&1/2&&&&&&&&& &  & \\
2T5 &&1/3&&&&1/6&&&&& &  & \\
2T6 &&1/6&&&&&1/3&&&& &  & \\
2O1 & 1/12&5/12&&&&&&&& &  & &\\
2O2 & 1/6&1/3&&&&&&&& &  & &\\\
2O3 & 1/3&1/6&&&&&&&& &  & &\\
2O4 & 5/12&1/12&&&&&&&& &  & &\\
2O5 & 1/3&1/6&&&&&&&& &  & &\\
2O6 & 1/6&1/3&&&&&&&& &  & &\\
2+ & &&1/36&2/36&3/36&2/36&1/36&3/36&2/36&1/36&&1/36&2/36\\
2-&1/4&&&&&&&&1/36&2/36&3/36&2/36&1/36
\end{array}
$
\caption{Transition matrix of the friend chain (first part).}  
\label{Fig:3}
\end{figure}

\begin{figure}
$
\begin{array}{c|ccccccccccccc}
&2T2 & 2T3 & 2T4 & 2T5 & 2T6 & 2O1 & 2O2 & 2O3 & 2O4 & 2O5 & 2O6 & 2+ & 2-\\
\hline
1+  &&&&&&&&&&& & 1/4 & 1/4\\
1- & &&&&& 3/36 & 2/36 &1/36&&1/36&2/36&&1/4\\
1T2 &  1/3 &&&&&& 1/6 &&&&&&\\
1T3 & &1/6 &&&&&& 1/3 &&&&&\\
1T4 & &  &&&&&&& 1/2 &&&&\\
1T5 &&&& 1/6 &&&&&& 1/3 &&&\\
1T6 &&&&&  1/3 &&&&&& 1/6 &&\\
1O1 &1/24 &&  &&1/24&&&&  &&&10/24&\\
1O2 &1/12 &1/12&  &&&&&& &&&1/3&\\
1O3 &1/24 &1/6& 1/8 &&&&&& &&&1/6&\\
1O4 &&1/12& 3/12 &1/12&&&&& &&&1/12&\\
1O5 &&& 3/24 &1/6&1/24&&&& &&&1/6&\\
1O6 &&& &1/12&1/12&&&& &&&1/3&\\

2T2  &1/3&&&&&&&&&& & 1/6 & \\
2T3  &&1/6&&&&&&&&& & 1/3 & \\
2T4  &&&&&&&&&&& & 1/2 & \\
2T5  &&&&1/6&&&&&&& & 1/3 & \\
2T6  &&&&&1/3&&&&&& & 1/6 & \\
2O1 &1/24&&&&1/24&1/4&1/12&&&& 1/12&  & \\
2O2 & 1/12&1/12&&&&3/24&1/6&1/24&& &  & &\\
2O3 & 1/24&1/6&3/24&&&&1/12&1/12&& &  & &\\
2O4 &&1/12&3/12&1/12&&&&1/24&&1/24& &   &\\
2O5 &&&3/24&1/6&1/24&&&&& 1/12&1/12   &\\
2O6 &&&&1/12&1/12&3/24&&&& 1/24&1/6  &\\
2+ &&&&&&&&&& & &1/4 & 1/4\\
2-&&&&&&&1/36&2/36&3/36&2/36&1/36&&1/4
\end{array}
$ 
\caption{Transition matrix of the friend chain (second part).}  
\label{Fig:4}
\end{figure}

\subsection{Numerics for the friend chain}
Let us investigate the mixing properties of the friend chain.
Standard theory quantifies ``closeness to stationarity after $n$ steps" via {\em variation distance} $d^*(n)$ or {\em separation distance} $s^*(n)$ from worst-case start, that is via
\begin{eqnarray*}
d^*(n) & := & \max_i \sfrac{1}{2} \sum_j   | p^n_{ij} - \pi_j |\\
s^*(n) & := & \max_{i,j} (1 - p^n_{ij}/\pi_j) 
\end{eqnarray*}
for the $n$-step transition matrix  $\bP^n$. Another measure of distance to stationarity is the $L^2$ or the $\chi^2$ distance. The $L^2$ distance between $P_i$ and the stationary distribution $\pi$ after $n$ steps is
\[ \| P^n(i, \cdot) - \pi \|_2 = \sqrt{\sum_j \frac{(p_{ij}^n - \pi_j)^2}{\pi_j}}\]
and 
\[ \| P^n - \pi \|_2 = \max _i    \| P^n(i, \cdot) - \pi \|_2     .\] 
Shown in Table \ref{Fig:5} are the numerical values for these distances.

\begin{table}[h]
$
\begin{array}{c|ccccccccc}
n &1&2&3&4&5&6&7&8&9\\
\hline
d^*(n) & 0.957 & 0.638 & 0.375 & 0.263 & 0.180 & 0.122 & 0.083 & 0.058 & 0.040\\
s^*(n) & 1 & 1 & 1 & 0.933 & 0.508 & 0.374 & 0.297 & 0.223 & 0.160\\
L^2(n) & 4.690 & 2.254 & 1.544 & 1.05 & 0.71 & 0.492 & 0.339 & 0.233 & 0.159
\end{array}
$ 
\caption{Measures of distance to stationarity for the friend chain, after $n$ games.}  
\label{Fig:5}
\end{table}

But these are not ``observable" quantities. More relevant to players is the mean number of games is which {\em friend} is on the same team, or on the opposing team,  as {\em ego}. This is a simple  calculation involving only matrix powers, and shown in Tables \ref{tab:meangames1} and \ref{tab:meangames2} are the mean number of games in which {\em friend} is on the opponent team and in which {\em friend} is on the same team.

\begin{table}[h!]
    \resizebox{0.9\textwidth}{!}{\begin{minipage}{\textwidth}
    \begin{tabular}{c|c c c c c c c c c c c c c c}
        Start & 1+ & 1- & 1T2 & 1T3 & 1T4 & 1T5 & 1T6 & 1O1 & 1O2 & 1O3 & 1O4 & 1O5 & 1O6\\ \hline 
        OT & 1.607 & 1.962 & 1.803 & 2.107 & 2.222 & 2.107 & 1.803 & 3.059 & 2.894 & 2.606 & 2.482 & 2.606 & 2.894 \\
        ST & 1.093 & 1.515 & 3.773 & 2.725 & 2.314 & 2.725 & 3.773 & 1.421 & 1.499 & 1.550 & 1.523 & 1.550 & 1.499 
    \end{tabular}
    \end{minipage}}
    \caption{Mean number of games (out of $8$) in which {\em friend} is on the opposite team OT (and the same team ST) as {\em ego}, who starts in the  first quadrant.}
    \label{tab:meangames1}
\end{table}
\begin{table}[h!]
    \resizebox{0.9\textwidth}{!}{\begin{minipage}{\textwidth}
    \begin{tabular}{c|c c c c c c c c c c c c c c}
        Start & 2T2 & 2T3 & 2T4 & 2T5 & 2T6 &  2O1 & 2O2 & 2O3 & 2O4 & 2O5 & 2O6 & 2+ & 2-\\  \hline 
        OT & 1.493 & 1.678 & 1.700 & 1.678 & 1.493 & 3.059 & 2.894 & 2.606 & 2.482 & 2.606 & 2.894 & 1.962 & 1.940 \\
        ST & 3.778 & 2.698 & 2.297 & 2.698 & 3.778 & 1.421 & 1.499 & 1.550 & 1.523 & 1.550 & 1.499 & 1.515 & 1.103
    \end{tabular}
    \end{minipage}}
    \caption{Mean number of games (out of $8$) in which {\em friend} is on the opposite team OT (and the same team ST) as {\em ego}, who starts in the second quadrant.}
    \label{tab:meangames2}
\end{table}

Notice that there is a symmetry property visible in Tables \ref{tab:meangames1} and \ref{tab:meangames2}, i.e. the values under $1T2$ to $1T6$, $1O2$ to $1O6$, $2T2$ to $2T6$, and $2O2$ to $2O6$ are exactly invariant under reversal. Moreover, the values under $1-$ and $2+$ are the same. This shows that if {\em friend} is playing with or against {\em ego}, the mean number of games in which {\em friend} is on the opposite team (or the same team) can be computed depending on the distance from {\em friend's} initial position to {\em ego} or {\em temporary ego's} initial position. In other words, we can reduce the number of states, and the states will depend on {\em friend's} shortest distance (counterclockwise or clockwise) from {\em ego}. 

A related question is the chance that you {\em never} play as an opponent (or as teammate) to your friend. Shown in Tables \ref{tab:notopp} and \ref{tab:notsame} are the numerical values for ``opponent" and ``teammate", respectively, omitting the cases where this is zero (initial opponent, or $1T4$).

\begin{table}[h!]
    \centering
    \begin{tabular}{c| c c c c c c c c}
        Start & 1+ & 1- &  1T2 & 1T3 & 1T5 & 1T6 \\  \hline 
        Probability & 0.098 & 0.057 & 0.141 & 0.026 & 0.026 & 0.141
    \end{tabular}
    \begin{tabular}{c| c c c c c c c c}
        Start & 2T2 & 2T3 & 2T4 & 2T5 & 2T6 & 2+ & 2-\\  \hline 
        Probability & 0.168 & 0.081 & 0.082 & 0.081 & 0.168 & 0.057 & 0.057
    \end{tabular}
    \caption{Probability that {\em friend} is never an opponent of {\em ego} over 8 games.}
    \label{tab:notopp}
\end{table}

\begin{table}[h!]
    \centering
    \begin{tabular}{c| c c c c c c c c}
        Start & 1+ & 1- & 1O1 & 1O2 & 1O3 & 1O4 & 1O5 & 1O6 \\  \hline 
        Probability & 0.403 & 0.292 & 0.344 & 0.317 & 0.271 & 0.251 & 0.271 & 0.317 
    \end{tabular}
    \begin{tabular}{c| c c c c c c c c}
        Start & 2O1 & 2O2 & 2O3 & 2O4 & 2O5 & 2O6 & 2+ & 2-\\  \hline 
        Probability & 0.344 & 0.317 & 0.271 & 0.251 & 0.271 & 0.317 & 0.292 & 0.393
    \end{tabular}
    \caption{Probability that {\em friend} is never a teammate of {\em ego} over 8 games.}
    \label{tab:notsame}
\end{table}

\paragraph{Comparison with random teams.}
There are several ways one could compare the effect of the ``mixing" scheme we study with the alternate scheme of randomly assigning players to team for every game.
For instance, under random mixing, if {\em friend} starts somewhere which is not the opposite team as {\em ego},
 then in each subsequent game  the probability that {\em friend} 
 is not on the opposite team as {\em ego} equals $17/23$.  So by independence,
  the probability that {\em friend} will never be on the opposite team as {\em ego} over eight games equals  $(17/23)^7 = 0.1205168...\approx 0.121$. 
 We see from Table \ref{tab:notopp} that this is less than the values under our scheme, if the starting position is only one position away from {\em ego} or the corresponding position of {\em ego} on the other half court. In the remaining cases, it is much greater. 
 
Similarly, if {\em friend} starts somewhere which is not the same team as {\em ego}, then (under random mixing) in each subsequent game  
  the probability that {\em friend} is not on the same team as {\em ego} equals $18/23$.
 So the probability that {\em friend} will never be on the same team as {\em ego} over eight games equals  $(18/23)^7 = 0.1798095...\approx 0.180$.
From Table \ref{tab:notsame}  this is always less than the values under our scheme.

\subsection{Monte Carlo Simulations}
More complicated ``observables" can most easily be addressed via Monte Carlo simulation of the process.
For instance
\begin{quote}
What is the probability that {\em ego} will encounter (as either teammate or opponent) {\em all} of the other $23$ players during an 8-game sequence?
\end{quote}
By Monte Carlo, the probability $\approx 0.595$ if  {\em ego} starts at the $1^\text{st}$ half court, or $\approx 0.675$ 
if {\em ego} starts in the $2^\text{nd}$ half court.
(Intuitively, these figures differ because in the former case one has more overlap between opponents in the first and second games.)
Over 10 games, these probabilities increase to $0.814$ and $0.857$. 
These simulations were done with 1 million trials, so we can take confidence intervals to be $\pm 0.001$.

\subsection{Solving Linear Systems}
\label{SLS}
Aside from computing powers of the transition matrix or using Monte Carlo simulations to answer numerical questions, 
a textbook method for calculation is via solving  linear systems of equations.  We will illustrate by two examples which can be done ``by hand" -- more complicated examples could be done numerically. 

\smallskip
\noindent
{\bf Example.} {\em For each of the various positions for players on the same team, what is the expected number of games until they are no longer on the same team?}

\smallskip
\noindent 
Recall the Figure \ref{Fig:26state} notation.
By symmetry, we can take the first player to be \emph{ego} in the first quadrant, and the other player to be one of $1T2, 1T3, 1T4$,
or  \emph{ego} in the second quadrant, and the other player to be one of $2T2, 2T3, 2T4$.
For one of those ``initial other player positions" $i$,
let $t_{i}$ be the expected number of games it takes until the two players are on different teams.
Consider $i = 1T2$. 
There are three equally likely outcomes of the first match.
Either neither of those two players migrates, in which case the configuration remains $1T2$.
Or both players migrate, in which case the configuration becomes $2T2$.
Or exactly one player migrates, in which case they are now on different teams. 
Repeating this analysis for each $i$ leads to
 the following linear system of equations:
\begin{align*}
    t_{1T2} &= 1 + \sfrac {1}{3} t_{1T2} + \sfrac {1}{3} t_{2T2} \\
    t_{2T2} &= 1 + \sfrac {1}{3} t_{1T2} + \sfrac {1}{3} t_{2T2} \\
    t_{1T3} &= 1 + \sfrac {1}{6} t_{1T3} +  \sfrac {1}{6}  t_{2T3} \\
    t_{2T3} &= 1 + \sfrac {1}{6} t_{1T3} +   \sfrac {1}{6}  t_{2T3} \\ 
    t_{1T4} &= 1 \\
    t_{2T4} &= 1.
\end{align*}
By symmetry we have $ t_{1T2} = t_{2T2},  t_{1T3} = t_{2T3} = t_{1T5} = t_{2T5},  t_{1T4} = t_{2T4}$.
Solving theses equations leads to
\begin{align*}
  t_{1T2} &= t_{2T2} = 3\\
  t_{1T3} &= t_{2T3} = t_{1T5} = t_{2T5} = \sfrac {3}{2}\\
  t_{1T4} &= t_{2T4} = 1.
\end{align*}
In fact this example could be solved without writing down equations, by first observing that in each case the number of games required has a Geometric($p$) distribution,
with $p = 1/3, 2/3,1$ according as the initial distance between players being $1, 2, 3$.


\smallskip
\noindent
{\bf Example.}
{\em Starting the friend Markov chain from one of the states $1 \pm$ or $2 \pm$, what is the probability that \emph{ego} and \emph{friend} will play on the same team before they play on opposing teams?}

\smallskip
\noindent
Here we do need to write down the equations.
Let $p_{i}$ be the probability that \emph{ego} and \emph{friend} will play on the same team before they play on opposing teams if \emph{friend} starts on state $i$, where $i$ is one of $1 \pm$ and $2 \pm$. 
While they are on opposite courts, there are $4$ equally likely possibilities for the players to migrate or not.
From the transition matrix of the friend chain, we obtain the following linear system:
\begin{align*}
    p_{1+} &= \sfrac{1}{4} p_{1+} +  \sfrac{1}{4} p_{1-} +  \sfrac{1}{4} p_{2+} +  \sfrac{1}{4} p_{2-}  \\
    p_{1-} &=  \sfrac{1}{4} +  \sfrac{1}{4} p_{1-} +  \sfrac{1}{4} p_{2-} \\
    p_{2+} &=  \sfrac{1}{4} +  \sfrac{1}{4} p_{2+} +  \sfrac{1}{4} p_{2-} \\
    p_{2-} &=  \sfrac{1}{4} p_{1+} +  \sfrac{1}{4} p_{2-} .
\end{align*}
The solution is
\[
    p_{1+} = \sfrac{3}{11} , \quad 
   p_{1-} = p_{2+} = \sfrac {4}{11}, \quad
    p_{2-} = \sfrac {1}{11}.
\]
Hence, if \emph{friend} starts on one of the states in $1+$, $1-$, $2+$, and $2-$, the probabilities that \emph{ego} and \emph{friend} will play on the same team before they play on opposing teams are $\frac{3}{11}$,$ \frac 4{11}$, $ \frac 4{11}$, and $\frac 1{11}$, respectively. 

\subsection{Suggestions for  course projects}
\label{sec:projects}

As mentioned in the introduction, one could use this topic as a source of projects to accompany a course in Markov chains with some emphasis on computation.
Here are some suggestions.

{\bf 1.} Repeat the analysis for the simpler model is which there are only 12 players and one court.
At the end of each match, the players in one front row swap positions with the players in the other front row.

{\bf 2.} Similar to the two examples in section \ref{SLS} above, 
calculate the expected amount of time until {\em ego} and {\em friend} are on the same team (or opposite teams), for each possible initial configuration.

{\bf 3.} (Suggested by a referee) 
Calculate the {\em fundamental matrix} of the friend chain.  The fundamental matrix determines the mean time to go from one given initial state to another given target state.  
Moreover, in the context of ``comparison with random teams" (the number of games you play with or against your friend in the Markov scheme, versus the number in the purely random model), the fundamental matrix determines the $n \to \infty$ limit of the expected difference between the two schemes.

{\bf 4.} A more challenging theory project is to improve our bounds on the mixing time.  The construction in section \ref{sec:irreducible}  implicitly gives a (very large) upper bound, and
the next section gives lower bounds, but these are surely far from optimal.

 \section{Mixing time for the big chain}
 \label{sec:mixing}
At a research level there has been extensive study of mixing times for many different  card-shuffling models, usually  in the asymptotic setting (as the size $n$ of card deck $\to \infty$). Our model is rather specific to the $n = 24$ case, so we have not sought to embed it into some family allowing large $n$.

One can get lower bounds on mixing time by considering specific functions of the chain, and the variation distances in Table \ref{Fig:5} for the friend chain are a lower bound for the distances in the big chain.
A first step \cite{me-RW} in studying some ``local moves" shuffles such as random adjacent  transpositions
was to obtain lower bounds by studying motion of initially adjacent cards.
In our model, the ``friends" maximal-start variation distance in Table \ref{Fig:5}  is indeed\footnote{Except for opening games.} from the case of initial adjacent players.

Recall that the progress of ego around the 4 quadrants is just the lazy cyclic walk, for which variation distance to stationarity in recorded in Table \ref{Fig:8}.

\begin{table}[h]
$
\begin{array}{c|ccccccccc}
n &1&2&3&4&5&6&7&8&9\\
\hline
d(n) & 0.5 & 0.25 & 0.25 & 0.125 & 0.125 & 0.0625 & 0.0625 & 0.0313 & 0.0313\\
\end{array}
$
\caption{Variation distance for the lazy cyclic walk.} \label{Fig:8}
\end{table}

These values are slightly less than the values in Table \ref{Fig:5} for the friend chain -- both are {\em priori} lower bounds for the variation distance for the big chain.

To find a better lower bound for the mixing time for the ``big" chain, we can combine the two aforementioned ideas and form a $52$- state (``big friend") chain, where {\em ego} can be in any of the 4 quadrants. In this new chain, we label the positions with respect to {\em ego} similar to the friend chain (see Figures \ref{Fig:26state} and \ref{Fig:52state}). We first record which quadrant {\em ego} is in (denoted as 1, 2, 3 or 4) and then whether friend is in the same team (denoted T) or the current opposing team (denoted O), or on the other court. If on the other court, we only note which quadrant, so denote by one of + or -. Finally, writing temporarily ego* for the opponent in the same position as ego, we indicate a friend's position as 1, 2, 3, 4, 5, or 6, counter-clockwise from ego or ego*. Its stationary distribution $\pi$ is induced from the uniform stationary distribution of the ``big friend" chain: probability $\frac6{92}$ for each of the states $1+, 1-, 2+$, $2-$, $3+$, $3-$, $4+$, and $4-$, and probability $\frac1{92}$ for each of the $44$ remaining states.

\setlength{\unitlength}{0.20in}
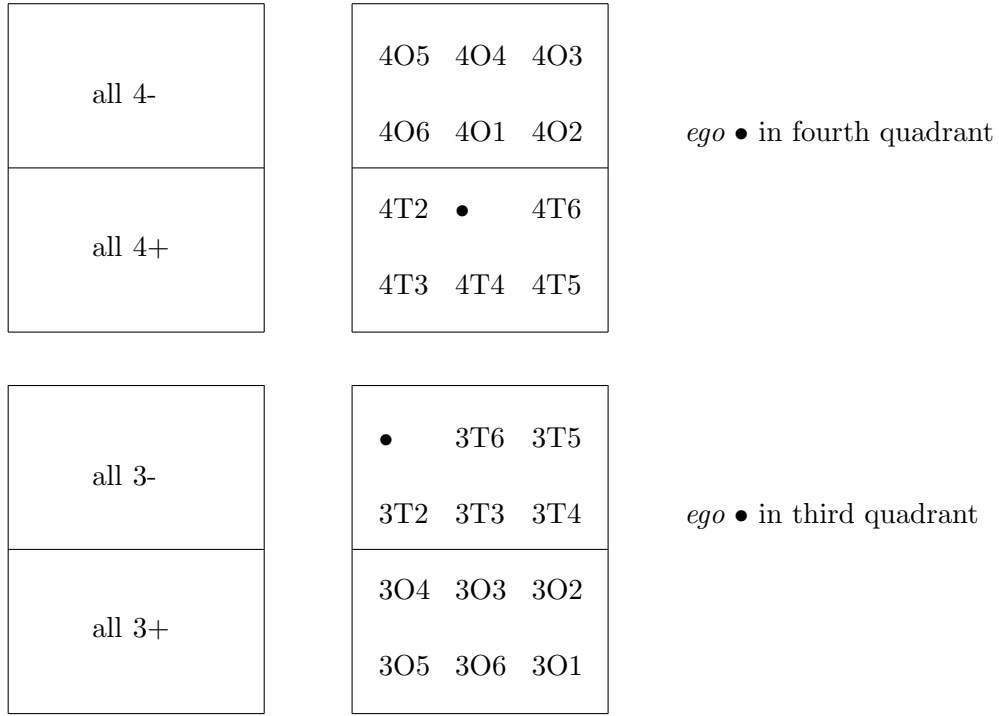
\begin{figure}[h!]
\begin{picture}(18,28)(0,-8)
\put(0.3,0){\line(1,0){6.7}}
\put(0.3,4.3){\line(1,0){6.7}}
\put(0.3,8.6){\line(1,0){6.7}}
\put(0.3,0){\line(0,1){8.6}}
\put(7,0){\line(0,1){8.6}}
\put(10,1){3O5}
\put(12,1){3O6}
\put(14,1){3O1}
\put(10,3){3O4}
\put(12,3){3O3}
\put(14,3){3O2}
\put(10,5){3T2}
\put(12,5){3T3}
\put(14,5){3T4}
\put(10,7){$\bullet$}
\put(12,7){3T6}
\put(14,7){3T5}
\put(18,5){{\em ego} $\bullet$ in third quadrant}
\put(9.3,0){\line(1,0){6.7}}
\put(9.3,4.3){\line(1,0){6.7}}
\put(9.3,8.6){\line(1,0){6.7}}
\put(9.3,0){\line(0,1){8.6}}
\put(16,0){\line(0,1){8.6}}
\put(2.5,2){all 3+}
\put(2.5,6){all 3-}
\put(0.3,10){\line(1,0){6.7}}
\put(0.3,14.3){\line(1,0){6.7}}
\put(0.3,18.6){\line(1,0){6.7}}
\put(0.3,10){\line(0,1){8.6}}
\put(7,10){\line(0,1){8.6}}
\put(10,11){4T3}
\put(12,11){4T4}
\put(14,11){4T5}
\put(10,13){4T2}
\put(12,13){$\bullet$}
\put(14,13){4T6}
\put(10,15){4O6}
\put(12,15){4O1}
\put(14,15){4O2}
\put(10,17){4O5}
\put(12,17){4O4}
\put(14,17){4O3}
\put(18,15){{\em ego}  $\bullet$  in fourth quadrant}
\put(9.3,10){\line(1,0){6.7}}
\put(9.3,14.3){\line(1,0){6.7}}
\put(9.3,18.6){\line(1,0){6.7}}
\put(9.3,10,0){\line(0,1){8.6}}
\put(16,10,0){\line(0,1){8.6}}

\put(2.5,12){all 4+}
\put(2.5,16){all 4-}

\end{picture}
\vspace*{-1.5in}
\caption{Positions of {\em friend} relative to {\em ego} $\bullet$.}  
\label{Fig:52state}
\end{figure}

The transition matrix for this ``big friend" chain is easy to obtain from the transition matrix of the ``friend" chain. Consider the transition matrix of the ``friend" chain in Figures \ref{Fig:3} and \ref{Fig:4}. We can view it as a matrix composed of smaller matrices as in Table \ref{tab:smallmatfriend}, where $T_{11}, T_{12}, T_{21},$ and $T_{22}$ are $13 \times 13$ matrices. The states used in the transition matrix for the  ``big friend" chain are arranged as follows: $1+$, $1-$, $1T2$ to $1T6$, $1O1$ to $1O6$, $2T2$ to $2T6$, $2O1$ to $2O6$, $2+$, $2-$, $3+$, $3-$, $3T2$ to $3T6$, $3O1$ to $3O6$, $4T2$ to $4T6$, $4O1$ to $4O6$, $4+$, and $4-$, and transition matrix for the  ``big friend" chain is shown in Table \ref{tab:smallmatbigfriend}.
\begin{table}[h!]
    \centering
    \begin{tabular}{c| c c}
        & 1 & 2 \\ \hline 
        1 & $T_{11}$ & $T_{12}$ \\ \hline 
        2 & $T_{21}$ & $T_{22}$
    \end{tabular}
    \caption{Transition matrix for the friend chain}
    \label{tab:smallmatfriend}
\end{table}

\begin{table}[h!]
    \centering
    \begin{tabular}{c| c c c c c}
        & 1 & 2 & 3 & 4 \\ \hline 
        1 & $T_{11}$ & $T_{12}$ & 0 & 0 \\ \hline 
        2 & $0$ & $T_{22}$ & $T_{21}$ & 0\\ \hline 
        3 & $0$ & 0 & $T_{11}$ & $T_{12}$ \\ \hline 
        4 & $T_{21}$ & 0 & 0 & $T_{22}$
    \end{tabular}
    \caption{Transition matrix for the big friend chain}
    \label{tab:smallmatbigfriend}
\end{table}

We can now investigate the mixing properties of the big friend chain. Shown in Table \ref{tab:l2distbig} are the distances from stationarity for the big chain.
\begin{table}[h!]
    \centering
    \begin{tabular}{c|c c c c c c c c c}
        n & 1 & 2 & 3 & 4 & 5 & 6 & 7 & 8 & 9\\ \hline 
        $d^*(n)$ & 0.978 & 0.713 & 0.520 & 0.340 & 0.242 & 0.168 & 0.125 & 0.085 & 0.058\\
        $s^*(n)$ & 1 & 1 & 1 &  1 & 0.827 & 0.681 & 0.461 & 0.391 & 0.272\\
        $L^2(n)$ & 6.708 & 2.977 & 1.868 & 1.228 & 0.827 & 0.563 & 0.387 & 0.266 & 0.183
    \end{tabular}
    \caption{Measures of distance to stationarity for the ``big friend" chain, after $n$ games}
    \label{tab:l2distbig}
\end{table}

Notice that the values we got for $d^*(n)$ in the ``big friend" chain are larger than that of the ``friend" chain, improving the lower bound of the ``big friend" chain.

\comment{
This, and other ideas for lower bounds could be more student projects.
For instance consider 3 adjacent players starting in the first quadrant (as in Figure \ref{Fig:26state}).
What is the probability that,
after exactly 5 games, 2 or 3 of these players are still on the same court (either half)?
Under the uniform distribution this event would have probability 1/2, by a symmetry argument.
In the actual chain, simulations show the probability is approximately 0.311.
This implies a  lower bound of approximately 0.189 for $d^*(5)$, slightly improving the bound in Table \ref{Fig:5}.
}

\section{Final remarks}
\label{sec:final}
We have mentioned analogies with card shuffling several times, because
 ``rotations" 
of team players correspond to a cut-shuffle of a 6-card deck.  
Our model is equivalent to a certain (not very easily implemented physically) random shuffle of a 24-card deck via first breaking into 4 sub-decks.
Persi Diaconis (personal communication) remarks that casinos and some fantasy games involve shuffling decks much larger than the usual 52-card deck,
and this is often done via some scheme involving  breaking into sub-decks, shuffling each in some way, and recombining in some way.  
Such schemes (thereby loosely analogous to our model) have generally not been studied in mathematical probability, an exception being the 
``casino shelf shuffling machines" studied by Diaconis-Fulman-Holmes \cite{shelf-shuffle}.

We have  interpreted the underlying question
\begin{quote} 
how effective is this scheme at mixing up the teams?
\end{quote}
in terms of mixing times, that is implicitly by comparison with the alternative of randomly assigning players to teams for each game.
An opposite alternative would be some analog of ``design of statistical experiments" schemes, deterministically assigning players to teams in each round  
in such a way that relative positions of two players were as uniformly spread as possible.  
At a practical level, our scheme is much easier and faster to implement than either alternative.
Moreover its implementation is {\em robust } to small variation in number of players, which is common in informal settings:  an extra player will rotate off court, or a 5-player team always has 3 players deemed front row.

\paragraph{Acknowledgements.} We thank an anonymous referee for helpful suggestions for exposition.


\end{document}